\documentclass[14pt]{article}
\usepackage[14pt]{extsizes}
\usepackage[utf8]{inputenc}
\usepackage{amsmath,amssymb,amsthm}
\usepackage{graphicx}

\usepackage{authblk}

\title{Asymptotics of Fundamental Solution of Cauchy Problem for Parabolic Equation with Small Parameter and Degeneration}
\author{Mark Rakhel}
\affil{National Research University Higher School of Economics, Moscow, Russia (marakhel@edu.hse.ru)}
\date{}

\begin{document}
  \maketitle

\begin{abstract}

In this paper, the method of constructing the asymptotics of the fundamental solution of the Cauchy problem for a degenerate linear parabolic equation with small diffusion is considered. Based on the results obtained in \cite{dn}, the study extends them over the case of a degenerate equation. As in \cite{dn}, the main technique that allows us to switch from pseudo-differential equations to partial differential equations is the non-oscillating WKB method. A distinctive feature of this work is a more detailed consideration on the characteristics of the Green's function in terms of symplectic geometry. The most significant intermediate result is presented as a theorem on the properties of the fundamental solution.

\end{abstract}

\newpage

\section{Introduction}

One of the powerful methods for constructive study of partial differential equations is the construction of asymptotic solutions. In the theory of hyperbolic equations, methods for constructing asymptotic solutions are distinguished, actively developed and applied, namely: the theory of the canonical Maslov operator, the theory of Fourier integral operators. These theories use objects and concepts of symplectic geometry: Lagrangian manifolds and symplectic forms in phase space. This theory is less developed for parabolic equations, a version of the theory of the Maslov canonical operator for parabolic equations is known, but in general this area needs to be studied. The fact is that the natural apparatus connecting hyperbolic equations with objects of symplectic geometry is integral transforms, in particular, the Fourier transform. For parabolic equations the use of this apparatus leads to the appearance of complex-valued functions that are not natural in the theory of parabolic equations. Nevertheless,  another approach can be used here. It is based on the representation of the Dirac delta function using the action of the function of the creation and annihilation operators on a Gaussian exponent in some abstract Hilbert space. This approach allows one to transfer the method developed in the theory of hyperbolic equations to parabolic equations. Similar problems in the case of non-degenerate diffusion were solved in the classical works of S. Varadan \cite{var}, S. Molchanov \cite{mol}, V. Maslov \cite{mas1} and also in the work of V. Maslov and V. Nazaikinskii \cite{masnaz} on the contemporary level relatively recently. The approach developed in these works essentially uses the non-degeneracy of the metric corresponding to the main symbol of the parabolic equation. Since this property is not present in particular problem, the technique mentioned above is applied.

\section{Delta function representations}

$\quad$ The key point in constructing asymptotics is the representation of the Dirac delta function via the action of the function of the creation operator on a Gaussian exponent, which allows us to make the conversion from the Cauchy problem for the fundamental solution to the Cauchy problem for the so-called symbol of fundamental solution. The formula itself has been introduced and proved in \cite{dn}.

\subsection{Delta function via Hermite polynomials}
$\quad$ Let us consider one of the definitions of the Hermite polynomials
\begin{equation}\label{ERM}
    H_n(x) = \biggm[ \biggm(\frac{\partial}{\partial z}\biggm)^2 e^{xz - \frac{z^2}{2}} \biggm] \biggm|_{z=0}
\end{equation}

As it is well known, these polynomials form an orthogonal basis in $L_2$ with weight function $e^{-\frac{x^2}{2}}$. Then, in particular, $\forall \psi(x) \in C_0^\infty:$

\begin{equation}
    \psi(x) = \frac{1}{\sqrt{\pi}}\sum_{n=0}^\infty\frac{1}{n!}\psi_n H_n(x) e^{-\frac{x^2}{2}}, 
    \nonumber
\end{equation}
where $\psi_n=\langle \psi(x),H_n(x)e^{-\frac{x^2}{2}}\rangle = \int_{-\infty}^{+\infty}\psi(x)H_n(x)e^{-\frac{x^2}{2}}dx$ - inner product. 
So we have the relation
\begin{equation*}
\begin{gathered}
    \psi(x) = \frac{1}{\sqrt{\pi}}\sum_{n=0}^{\infty}\frac{1}{n!}\langle \psi(\xi),H_n(\xi)e^{-\frac{\xi^2}{2}} \rangle H_n(x)e^{-\frac{x^2}{2}} = \\ = \langle \psi(\xi),\frac{1}{\sqrt{\pi}} \sum_{n=0}^{\infty}\frac{1}{n!}H_n(x) H_n(\xi)e^{-\frac{x^2 + \xi^2}{2}} \rangle
\end{gathered}
\end{equation*}
By the definition, $\langle \delta(\xi - x),\psi(\xi) \rangle = \psi(x)$. This implies
\begin{equation*}
    \delta(\xi - x) = \frac{1}{\sqrt{\pi}} \sum_{n=0}^{\infty}\frac{1}{n!}H_n(x) H_n(\xi)e^{-\frac{x^2 + \xi^2}{2}}
    \nonumber
\end{equation*}
After considering $x = 0$ and redefining variables, relation looks as follows 
\begin{equation}
    \delta(x) = \frac{1}{\sqrt{\pi}} \sum_{n=0}^{\infty}\frac{1}{n!}H_n(0) H_n(x)e^{-\frac{x^2}{2}}
\end{equation}
Knowing that $e^A = \sum_{n=0}^{\infty}\frac{1}{n!} A^n$, we can introduce a small parameter $h$ and obtain
\begin{equation}
    \delta(x) = \frac{1}{\sqrt{\pi h}}e^{h\frac{\partial^2}{\partial y \partial z}} e^{-\frac{(x-z)^2 + y^2}{2h}}\biggm|_{z = y = 0} \stackrel{def}{=}\frac{1}{\sqrt{\pi h}}\mathcal{L}_h e^{-\frac{(x-z)^2 + y^2}{2h}}\biggm|_{z = y = 0},
\end{equation}
where $\mathcal{L}_h = e^{h\frac{\partial^2}{\partial y \partial z}}$.
\subsection{Delta function via creation-annihilation operators}
$\quad$ Creation and annihilation operators $a$ and $a^+$ are satisfying the commutation relation $[a,a^+] = h$.
For example, 
\begin{equation}\label{CA}
    \begin{cases}
        a = \zeta + \frac{h}{2}\frac{\partial}{\partial \zeta} \\
        a^+ = \zeta - \frac{h}{2}\frac{\partial}{\partial \zeta}  
\end{cases}
\end{equation}
We can define functions of these operators arranged in particular order
\begin{equation*}
    \begin{cases}
        f(\stackrel{2}{a^+},\stackrel{1}{a}) = \int{\tilde{f}(\alpha_1,\alpha_2)e^{-i\alpha_1a^+}e^{-i\alpha_2a}d\alpha_1d\alpha_2} \\
        g(\stackrel{1}{a^+},\stackrel{2}{a}) = \int{\tilde{g}(\alpha_1,\alpha_2)e^{-i\alpha_2a}e^{-i\alpha_1a^+}d\alpha_1d\alpha_2}
    \end{cases}
\end{equation*}
where $\tilde{f}(\alpha_1,\alpha_2)$ and $\tilde{g}(\alpha_1,\alpha_2)$ - are the Fourier transforms of the functions $f(z,y)$ and $g(z,y)$ 
\newline
We use \cite{mas}
\begin{equation}
    F(\stackrel{2}{a},\stackrel{1}{a^+}) = 	(\mathcal{L}_h F)(\stackrel{1}{a},\stackrel{2}{a^+})
    \nonumber
\end{equation}
Denote $F(z,y) = e^{\frac{-(x-z)^2-y^2}{2h}}$. Then
\begin{equation}
    \mathcal{L}_h e^{-\frac{(x-\stackrel{1}{a})^2}{2h} - \frac{(\stackrel{2}{a^+})^2}{2h}} = e^{-\frac{(x-\stackrel{2}{a})^2}{2h} - \frac{(\stackrel{1}{a^+})^2}{2h}} = e^{-\frac{(x-a)^2}{2h}}e^{-\frac{(a^+)^2}{2h}}
\end{equation}
Let us consider $\chi \in \ker(a), \quad \chi_+ \in \ker{(a^+)^*}$. Then $a\chi = 0, \quad (a^+)^*\chi_+ = 0$. We found out that $\delta(x) = \frac{1}{\sqrt{\pi h}}(\mathcal{L}_h F)(0,0)$. On the other hand
\begin{equation}\label{IP}
\begin{gathered}
    (\mathcal{L}_h F)(0,0) = \langle \chi_+,((\mathcal{L}_h F)(\stackrel{1}{a},\stackrel{2}{a^+})\chi\rangle \langle \chi_+,\chi\rangle^{-1} = \\ =
    \langle \chi_+,F(\stackrel{2}{a},\stackrel{1}{a^+})\chi\rangle \langle \chi_+,\chi\rangle^{-1}
\end{gathered}
\end{equation}

The null vectors of operators (\ref{CA}) are $\chi = \chi_+ = e^{-\frac{\zeta^2}{h}}, \quad \langle \chi_+,\chi\rangle^{-1} = \sqrt{\frac{2}{\pi h}}$. Substituting those expressions in (\ref{IP}), we get the key expression for the Dirac delta. 
\begin{equation}\label{DELTA}
    \delta(x) = \frac{\sqrt{2}}{\pi h}\langle e^{-\frac{(x-a^+)^2}{2h}}e^{-\frac{\zeta^2}{h}}, e^{-\frac{(a^+)^2}{2h}}e^{-\frac{\zeta^2}{h}}  \rangle
\end{equation}

\textbf{Note.} 
Instead of using (\ref{DELTA}), we can rely on the following lemma, which was initially presented and proven in \cite{dn}. 

\textbf{Lemma}. Any function $\phi(x,h)$ uniformly bounded for each fixed $h$ and satisfying the  Lipschitz condition $|\phi(x_1,h)-\phi(x_2,h)| \leq k|x_1-x_2|$ satisfies the relation

\begin{equation}
    \phi(0,h) = \frac{1}{\sqrt{2\pi h}} \lim_{\beta \to 1-0} \int_{-\infty}^{\infty} \phi(x,h) \langle e^{-\beta \frac{(x-a^+)^2}{2h}}e^{-\frac{\zeta^2}{h}}, \delta(\zeta)  \rangle dx
\end{equation}

\section{Algorithm}

$\quad$The Green's function $G(x,\xi,t,h)$ is defined as the solution to the Cauchy problem
\begin{equation}\label{EQ1}
    \begin{cases}
        -h\frac{\partial G}{\partial t} + \hat{H}(\stackrel{2}{x},\stackrel{1}{-h\frac{\partial}{\partial x}})G = 0 \\
        G|_{t=0} = \delta(x-\xi)
    \end{cases}    
\end{equation}
where $H(x,p)$ is a function analytic in $p$ and taking real values whenever the arguments are real.\\

From the lemma follows that if $V_{\beta}(x,\xi,y,t,h)$ is a solution to the Cauchy problem
\begin{equation}\label{EQ2}
    \begin{cases}-h\frac{\partial V_{\beta}}{\partial t} + H(\stackrel{2}{x},\stackrel{1}{-h\frac{\partial}{\partial x}})V_{\beta} = 0 \\
    V_{\beta}|_{t=0} = e^{-\beta\frac{(x-\xi-y)^2}{2h}}
    \end{cases}    
\end{equation}
then the general form of the fundamental solution can be represented by the relation
\begin{equation}
    G(x,\xi,t,h) = \frac{1}{\sqrt{2\pi h}}\lim_{\beta \to 1-0} V_{\beta}(x,\xi,a^+,t,h)e^{-\frac{\zeta^2}{h}}\biggm|_{\zeta = 0}
\end{equation}
Function $V(x,t,y,h)$ is called by a symbol of the fundamental solution and its argument $y$ corresponds to the creation operator $a^+$.

It is proposed to search for the symbol using the non-oscillating WKB method, that is, to look for asymptotic approximation $V_{\beta,N}(x,\xi,y,t,h)$
\begin{equation}
    V_{\beta,N} = e^{-\frac{\Phi_{\beta}}{h}}(\varphi_0 + h\varphi_1 + ... + O(h^N))
    \nonumber
\end{equation}
In this paper, we restrict ourselves to the first term
\begin{equation}\label{FIRST}
    V_{\beta} = e^{-\frac{\Phi_{\beta}}{h}}(\phi_0 + O(h))
\end{equation}
Hence, for function $\Phi_{\beta}(x,\xi,y,t)$ we get the Hamilton-Jacobi equation
\begin{equation}\label{HJ}
    \begin{cases}
        \frac{\partial \Phi_{\beta}}{\partial t} + H(x,\frac{\partial\Phi_{\beta}}{\partial x}) = 0 \\
        \Phi_{\beta}|_{t=0}=\beta \frac{(x-\xi-y)^2}{2}
    \end{cases}
\end{equation}
and for $\varphi_0(x,\xi,y,t)$ - Transport equation
\begin{equation}\label{TR}
    \begin{cases}
        \Pi \varphi = 0 \\
        \varphi_0|_{t=0} = P(x)
    \end{cases}
\end{equation}
where $\Pi = \frac{\partial}{\partial t} + H_p \frac{\partial}{\partial x} + \frac{1}{2}H_{pp}(\Phi_{\beta})_{xx}$, 
and function $P(x) \in C_0^\infty$ looks as following $$P(x) = \begin{cases}
            1, |x| < \delta_1 \\
            0, |x| \geq \delta 
        \end{cases}$$ if $\delta < \delta_1$. This means that we restrict ourselves to considering solutions in the neighborhood of $ x = 0 $ and initial conditions with support from the same neighborhood.
The characteristic system of (\ref{HJ}) is the Hamilton system
\begin{equation*}
    \begin{cases}
        \dot{x} = H_p, \quad x|_{t=0} = x_0 \\
        \dot{p} = -H_x, \quad p|_{t=0} = \beta(x_0 - \xi -y)
    \end{cases}
\end{equation*}

and the solution to the Hamilton-Jacobi equation can be found with the integration of form $\omega = pdx$ and projection $x_0(x,\xi,y,t)$ as long as $J_0 = |\frac{\partial x}{\partial x_0}(x_0,\xi,y,t)| \ne 0$ 

The expression for the solution to the Cauchy problem (\ref{TR}) is well known
\begin{equation}
    \varphi_0 = \frac{1}{\sqrt{J_0}}e^{\frac{1}{2}\int_0^t{H_{x p}d\tau}}
\end{equation}

The action of the function of the creation operator on the exponent before taking the limit, namely $V(x,\xi,a^+,t,h)e^{-\frac{\zeta^2}{h}}$ can be also figured out by dint of the non-oscillating WKB method. Eventually, the explicit formula for the fundamental solution is as follows

\begin{equation}\label{FUND}
    G(x,\xi,t,h) =  \frac{1}{\sqrt{2\pi h}}\lim_{\beta \to 1-0}{\frac{e^{\frac{1}{h}(\frac{\hat{y}^2}{2} - \Phi_{\beta}(x,\xi, \hat{y},t))}}{\sqrt{J_0}\sqrt{1 - \frac{\partial^2 \Phi_{\beta}}{\partial y^2}(x,\xi, \hat{y},t)}}}e^{\frac{1}{2}\int_0^t{H_{x p}d\tau}}(1+O(h))
\end{equation}
where $\hat{y} = \hat{y}(x,t,\xi)$ is the implicit function defined as a solution to the equation $y = \frac{\partial \Phi}{\partial y}(x,t,\xi)$.

\section{New geometric method}

$\quad$ A more detailed consideration of the construction described above indicates to us the presence of a new symplectic structure, an analogue of which arises in the theory of Fourier integral operators, namely:

1) Extended four-dimensional phase space with variables $(x,y,p_x,p_y)$ and Hamilton system
\begin{equation}
    \begin{cases}
        \dot{x}=H_p(x,p_x), \quad x|_{t=0}=x_0 \\
        \dot{p_x}=-H_x(x,p_x), \quad p_x|_{t=0}=\beta(x_0-y-\xi) \\
        \dot{y}=0, \quad y|_{t=0}=y\\
        \dot{p_y}=0, \quad p_y|_{t=0} = -\beta(x_0-y-\xi) 
    \end{cases}
\end{equation}
wherein $y$ and $p_y$ correspond to the creation operator. \newline
2) New fundamental form
\begin{equation}
    \omega_y \in T^*(\mathbb{R}_x^n) \times T^*(\mathbb{R}_y^n), \quad
    \omega_y = p_xdx - p_ydy
\end{equation}
3) Function $V$ appears in the variables $(x,t)$ after the projection
\begin{equation}
    \pi : \Lambda_t^{2n} \to (\mathbb{R}_x^n,y=p_y)
\end{equation}
which comes to be the product of the configuration space and the diagonal of the complement $\mathbb{R}_{y,p_y}^n$. \newline
As mentioned above, phase $\Phi$ is found by integrating the form and applying the inverse mapping
\begin{equation}
    \Phi_{\beta} = (\pi^{-1})^*\int \omega_y,
\end{equation}
where $(\pi^{-1})^*$ is the mapping induced by $\pi^{-1}$. \newline
It is easy to see that 
\begin{equation*}
\begin{gathered}
    \int{\omega}_y = \Phi_{\beta}|_{t=0} - \int_0^y{p_y dy} + \int_0^t{{(pH_p-H) d\tau}} = \\ = \frac{\beta}{2}(x_0 - \xi)^2 + \int_0^t{(pH_p-H) d\tau}    
\end{gathered}
\end{equation*}
If $(pH_p-H) \geq 0$, then for each $x_0$ and $y$
\begin{equation*}
    \Phi_{\beta}  \geq 0
\end{equation*}

So, according to the new notation the expression for the fundamental solution takes the form
\begin{equation}\label{NEWFUND}
     G(x,t,\xi,h) = \frac{1}{\sqrt{2\pi h}}\lim_{\beta \to 1-0}{\frac{e^{-\frac{1}{h}(\pi^{-1})^*\int \omega_y}}{\sqrt{J}}}e^{\frac{1}{2}\int_0^t{H_{xp}d\tau}}(1+O(h))
\end{equation}
where $J = |\frac{\partial(x,y-p_y)}{\partial(x_0,y)}|$

\section{Example: Heat equation}

$\quad$ Let us consider as an example the one-dimensional heat equation with constant coefficients and a small parameter $h$
\begin{equation}\label{HEATEXAMPLE}
    \frac{\partial u}{\partial t} - h\frac{\partial^2 u}{\partial x^2} = 0
\end{equation}
the fundamental solution of which is well known. Let us find it first using a sequential algorithm and then with the geometric method. \vspace{2mm} \\
The Green's function for (\ref{HEATEXAMPLE}) is the function $G(x,\xi,t,h)$, satisfying 
\begin{equation}
    \begin{cases}
        \frac{\partial G}{\partial t} - h\frac{\partial^2 G}{\partial x^2} = 0 \\
        G|_{t=0} = \delta(x-\xi)
    \end{cases}
\end{equation}
The problem for the symbol $V(x,\xi,y,t,h)$
\begin{equation}\label{SYMBHEAT}
    \begin{cases}
        \frac{\partial V}{\partial t} - h\frac{\partial^2 V}{\partial x^2} = 0 \\
        V|_{t=0} = e^{-\frac{(x-\xi - y)^2}{2h}}
    \end{cases}
\end{equation}
Non-oscillating WKB method leads to the Hamilton-Jacobi and the Transport equations for functions $\Phi(x,t,y)$ and $\varphi(x,t,y)$
\begin{equation}\label{PHIHEAT1}
    \begin{cases}
        \frac{\partial \Phi}{\partial t} + (\frac{\partial \Phi}{\partial x})^2 = 0 \\
        \Phi|_{t=0} = \frac{(x-\xi - y)^2}{2}
    \end{cases}
\end{equation}
The Hamilton system is 
\begin{equation*}
    \begin{cases}\label{HAM}
        \dot{x}=2p, \quad x|_{t=0}=x_0 \\
        \dot{p}=0, \quad p|_{t=0}=x_0-y-\xi
    \end{cases}
\end{equation*}
and the solution to the problem (\ref{PHIHEAT1}) is as follows
\begin{equation*}
    \Phi = \int{pdx}\big|_{x_0=x_0(x,\xi,y,t)} = \frac{(x-\xi-y)^2}{2(1+2t)}
\end{equation*}
A solution to the Transport equation $\varphi_0(x,\xi,y,t)$ in this case is
\begin{equation*}
    \varphi_0 = \frac{1}{\sqrt{\frac{\partial x}{\partial x_0}}} = \frac{1}{\sqrt{1 + 2 t}}
\end{equation*}
Now we can write down the explicit formula for $V$ 
\begin{equation}\label{SYMBHEAT1}
    V =  \frac{1}{\sqrt{1 + 2 t}}e^{-\frac{(x-\xi-y)^2}{2(1+2t)h}}
\end{equation}
Then in the expression for the symbol (\ref{SYMBHEAT1}), instead of $ y $, you need to substitute the creation operator $a^+ = \zeta - \frac{h}{2}\frac{\partial}{\partial \zeta}$ and calculate the action of the function from the operator on the Gaussian exponent $V(x,\xi,a^+,t,h)e^{-\frac{\zeta^2}{h}}|_{\zeta = 0}$.  This action is calculated through finding the function $ W (x, \ xi, t, h, \alpha, \zeta) $, which is a solution to the problem
\begin{equation}\label{HELPW}
    \begin{cases}
        \frac{\partial W}{\partial \alpha} - V(x,\xi,a^+,t,h)W = 0 \\
        W|_{\alpha = 0} = e^{-\frac{\zeta^2}{h}}
    \end{cases}
\end{equation}
Then $W(x,\xi,t,h,1,0) = V(x,\xi,a^+,t,h)e^{-\frac{\zeta^2}{h}}|_{\zeta = 0}$. The problem (\ref{HELPW}) is again solved using the non-oscillating WKB method. So, repeating the steps described above we obtain the required expression for the Green's function of the heat equation
\begin{equation}\label{HEATGREEN}
    G(x,\xi,t,h) = \frac{1}{\sqrt{4\pi h t}}e^{-\frac{(x-\xi)^2}{4 t h}}
\end{equation}

Now let's find a fundamental solution using the "geometric" method. In this case, instead of the system (\ref{HAM}), we write the extended system of characteristics
\begin{equation*}
    \begin{cases}
        \dot{x}=2p, \quad x|_{t=0}=x_0 \\
        \dot{p_x}=0, \quad p_x|_{t=0}=\beta(x_0-y-\xi) \\
        \dot{y}=0, \quad y|_{t=0}=y\\
        \dot{p_y}=0, \quad p_y|_{t=0} = -\beta(x_0-y-\xi) 
    \end{cases}
\end{equation*}
From this system, we now need to express $ x_0 (x, \xi, t) $ and $ y (x, \xi, t) $. The expressions for them are as follows
\begin{equation*}
\begin{gathered}
    x_0 = \frac{-x(1-\beta) - 2\beta\xi t}{1 - \beta + 2\beta t} \\
    y = \frac{-\beta(x-\xi)}{1-\beta + 2\beta t}
\end{gathered}
\end{equation*}
and the Jacobian of the extended system is
\begin{equation*}
    J = \frac{\partial(x,y-p_y)}{\partial(x_0,y)} = 1 - \beta + 2\beta t
\end{equation*}
Finally, by integrating the form $ \omega_y $ and substituting the required values into the formula (\ref{NEWFUND}), we obtain an expression for the Green's function of the heat equation (\ref{HEATGREEN}).

\section{Construction of the asymptotics of the fundamental solution for a linear degenerate parabolic equation}

$\quad$ Let us consider in details the case $ H (x, p) = x ^ 2p ^ 2 $, i.e., we seek a fundamental solution to the problem
\begin{equation}\label{EX}
    \begin{cases}
        \frac{\partial u}{\partial t} - h x^2 \frac{\partial^2 u}{\partial x^2} = 0 \\
        u|_{t=0} = u_0 (x)
    \end{cases}    
\end{equation}
namely, the function $ G (x, t, \xi, h) $, which is a solution to the Cauchy problem
\begin{equation}\label{FUNDDEG}
    \begin{cases}
        \frac{\partial G}{\partial t} - h x^2 \frac{\partial^2 G}{\partial x^2} = 0 \\
        G|_{t=0} = \delta(x-\xi)  
    \end{cases}    
\end{equation}
The solution to the original problem (\ref{EX}) is the function $ u (x, t) $, which is defined as the convolution of the Green's function with the initial data
\begin{equation*}
    u(x,t,h) = \int{G(x,\xi,t,h)u_0(\xi)d\xi}
\end{equation*}

\textbf{Note.}
Let $ \xi = 0 $. Then in the problem (\ref{FUNDDEG}) you can separate the variables, namely, look for a solution in the form $ G = f (t) \delta (x) $. Solving the ordinary differential equation, we get
\begin{equation*}
     G = e^{2ht}\delta(x) = \delta(x)(1+O(h))
\end{equation*}
This fact shows the absence of smoothing, in contrast to ordinary parabolic equations (for example, the heat equation).

The Cauchy problem for the symbol has the form
\begin{equation}
    \begin{cases}
        \frac{\partial V_{\beta}}{\partial t} - h x^2 \frac{\partial^2 V_{\beta}}{\partial x^2} = 0 \\
        V_{\beta}|_{t=0} = e^{-\beta\frac{(x-\xi - y)^2}{2h}}
    \end{cases}
\end{equation}
and after applying the WKB method, we obtain the problem for finding the phase $ \Phi _ {\beta} (x, \xi, y, t) $
\begin{equation}\label{HJE}
    \begin{cases}
        \frac{\partial \Phi_{\beta}}{\partial t} + x^2 (\frac{\partial \Phi_{\beta}}{\partial x})^2 = 0 \\
        \Phi_{\beta}|_{t=0}=\beta \frac{(x-\xi-y)^2}{2}
    \end{cases}
\end{equation}
and the amplitude $\phi_0(x,\xi,y,t)$
\begin{equation}
    \begin{cases}
        \frac{\partial \phi_0}{\partial t} + 2 (\Phi_{\beta})_x \frac{\partial \phi_0}{\partial x} + (\Phi_{\beta})_{xx}\phi_0 = 0 \\
        \phi_0|_{t=0} = P(x)
    \end{cases}
\end{equation}
Hamilton characteristic system for (\ref{HJE})
\begin{equation}
    \begin{cases}
        \dot{x} = H_p = 2x^2p, \quad x|_{t=0} = x_0 \\
        \dot{p} = -H_x = -2xp^2, \quad p|_{t=0} = \beta(x_0-\xi-y)
    \end{cases}
\end{equation}
From such a system we find $x = x(x_0,\xi,y,t)$ and $p=p(x_0,\xi,y,t)$:
\begin{equation*}
    \begin{cases}
        x = x_0 e^{2\beta x_0(x_0-\xi-y)t} \\
        p = \beta(x_0-\xi-y)e^{-2\beta x_0(x_0-\xi-y)t}
    \end{cases}
\end{equation*}
The Jacobian of the mapping $ J_0 = \frac {\partial x} {\partial x_0} (x_0, \xi, y, t) $ is here
\begin{equation*}
    J_0 = e^{2\beta x_0(x_0-\xi-y)t}(1 + 2\beta t x_0(2x_0-\xi-y))
\end{equation*}
If we do not introduce the parameter $ \beta $ (or immediately set it to one), such a Jacobian will take zero values for $ y = const $, and the Lagrangian manifold is S-shaped at some points. This is clearly seen at $t \approx 2/3$. \vspace{5mm}
\begin{figure}[h]
    \includegraphics[height = 8cm]{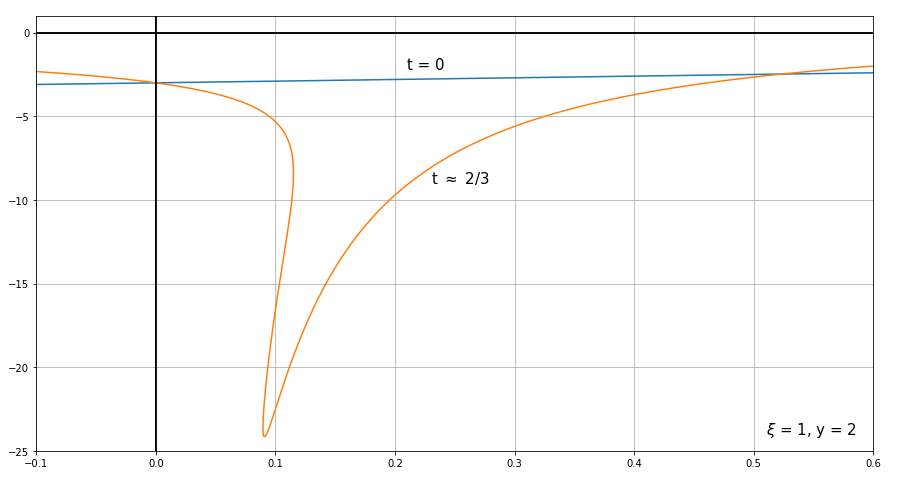}
\end{figure}
\vspace{1mm} \\

In terms of the new geometric method, a similar picture is the graph of the projection of the section of the Lagrangian manifold $ \lambda \in \mathbb {R} ^ 4 $ by the plane $ y = const $ onto the space $ (x, p_x) $. But in the very four-dimensional space $ (x, y, p_x, p_y) $ for sufficiently small $ t $ such pictures do not appear.

\textbf{Lemma.}
The function $ G _ {\beta} $ satisfies the following relations (in the case of distributions)
\begin{equation}
    \lim_{\beta \to 1-0}G_{\beta}(0,\xi,t,h) = \delta(\xi)
\end{equation}
and
\begin{equation}\label{DELTA111}
    \lim_{\beta \to 1-0} (\psi(x),G_{\beta}(x,0,t,h)) = \psi(0)(1+O(h)),
\end{equation}
$\psi$ - test function. 
\begin{proof}
    To prove the first relation we use the fact that if $ x = 0 $, then $ x_0 = 0 $. The function $ G _ {\beta} (0, \xi, t, h) $ is then written out explicitly
    \begin{equation*}
        G_{\beta} = \frac{e^{-\frac{\beta \xi^2}{2(1-\beta)h}}}{\sqrt{2\pi h}{\sqrt{1-\beta}}}
    \end{equation*}
    and it equals $\delta(\xi)$. \\
    The second relation can be proved using the following fact
    \begin{equation*}
        \int_{-\infty}^{\infty} G_{\beta}(x_0(x,\xi,t),\xi,t,h) dx = \int_{-\infty}^{\infty} G_\beta(x_0,\xi,t,h) \frac{dx}{dx_0} dx_0
    \end{equation*}
    The convergence of the integral is then easily verified, and again we get a delta function.
\end{proof}

\subsection{Existence and properties of the solution}

$\quad$ As noted in the relevant section, the improved algorithm introduces a special symplectic structure, including the projection $ \ pi $, which is the composition of the projection onto the $ x $ axis with the projection onto the diagonal $ y = p_y $. Thus, in the course of solving, it becomes necessary to express $ x_0 (x, \xi, \beta, t) $ and $ \hat {y} (x, \xi, \beta, t) $ from the system
\begin{equation*}
    \begin{cases}
        x(x_0,\xi,y,\beta,t) = x \\
        y = p_y(x_0,\xi,y,\beta,t)
    \end{cases}
\end{equation*}
For the case $ H (x, p) = x ^ 2p ^ 2 $ this system looks as follows
\begin{equation}\label{FUND123}
    \begin{cases}
        x = x_0 e^{2\beta x_0(x_0-\xi-y)t} \\
        y = -\beta(x_0-\xi-y)
    \end{cases}
\end{equation}
Expressing $ y (x_0, \xi, \beta, t) = - \frac {\beta (x_0- \xi)} {1- \beta} $ from the second equation of the system (\ref{FUND123}) and substituting into the first, we get an expression for $ x (x_0, \xi, \beta, t) $
\begin{equation}\label{EXP123}
    x = x_0 e^{\frac{2\beta x_0(x_0-\xi)t}{1-\beta}}
\end{equation}
$x_0(x,\xi,\beta,t)$, as noted above, cannot be expressed explicitly. But we can consider the Jacobian of the map by substituting $y = y(x0,\xi,\beta,t)$
\begin{equation}\label{JAC}
    J_0(x_0,\xi,\beta,t) = \frac{\partial x(x_0,\xi,\beta,t)}{\partial x_0} = e^{\frac{2\beta x_0(x_0-\xi)t}{1-\beta}}\big(1 + \frac{2\beta t}{1-\beta}(2x_0^2-x_0 \xi)\big)
\end{equation}

For $0 < \beta < 1$ for any bounded $x$ and $\xi$ this Jacobian is nowhere equal to zero for sufficiently small $ t $, which means that the solution will be a smooth function.

For $ \beta \to 1-0 $, the expression (\ref{EXP123}) is defined only for two values of $ x_0 $: $ x_0 = 0 $ and $ x_0 = \xi $, and for each of these values the Jacobian is defined and distinct from zero, which means that the solution also exists as a smooth function.

It is worth noting that the formula for the fundamental solution contains the Jacobian of the extended system $ J = \frac {\partial (x, y-p_y)} {\partial (x_0, y)} $, that is, (\ref{JAC}) is multiplied by $ \frac {\partial (y-p_y)} {\partial y} = 1- \beta $, and for $ \beta \to 1-0 $, the denominator no longer has an infinitely large factor.

Since $ x = 0 \Leftrightarrow x_0 = 0 $, this point is the support of the delta function, which remains at the point $ x = 0 $ for any $ t> 0 $. If $ x \not = 0 $, then any solution $ x ^ * $ of the equation of characteristics lies on the trajectory outgoing from the point $ x_0 = \xi $, which is the support of the initial data.

\subsection{asymptotics of the logarithmic limit of the fundamental solution as $t \to +0$}

$\quad$ To construct the asymptotics of the limit of the fundamental solution as $ t \to +0$, we need to introduce the following lemma \cite{df}

\textbf{Lemma.}
Suppose that there exists a smooth solution to the problem
\begin{equation}
    \begin{cases}
        \frac{\partial \Phi}{\partial t} + H(x,\frac{\partial\Phi}{\partial x}) = 0 \\
        \Phi|_{t=0}= \frac{(x-\xi-y)^2}{2}
    \end{cases}
\end{equation}
for $0 \leq t \leq T$, i.e. $J = \frac{\partial x(x_0,\xi,y,t)}{\partial x_0} \ne 0$. Then for $t>0$
\begin{equation}
     G(x,t,\xi,h) = \frac{1}{\sqrt{2\pi h}}{\frac{e^{-\frac{S(x,t,\xi)}{h}}}{\sqrt{J_0}}}e^{\frac{1}{2}\int_0^t{H_{xp}d\tau}}(1+O(h))
\end{equation}
where the function $S(x, t, \xi)$ is determined by the expression
\begin{equation}
     S(x,t,\xi) = \int_0^t{(pH_p(x,p,\tau) - H(x,p,\tau))d\tau}
\end{equation}
and $x = x(\xi, p_0, t)$, $p = p(\xi, p_0, t)$ is the solution of the Hamilton system
\begin{equation}
    \begin{cases}
        \frac{dx}{d\tau} = H_p(x,p,\tau), \quad x|_{\tau = 0} = \xi \\
        \frac{dp}{d\tau} = -H_x(x,p,\tau), \quad p|_{\tau = 0} = p_o \in \mathbb{R}
    \end{cases}
\end{equation}
and $J_0 = \frac{\partial x(x_0,\xi,p_o)}{\partial p_0} \ne 0$.

Next, since we want to investigate the behavior of the function $S(x,t,\xi)$ as $t \to +0$, we represent $t$ in the form $t = \nu t'$, where $t'$ is a fixed number and $\nu \to +0$. Changing the variable $t = \nu\sigma$, from the lemma we receive
\begin{equation}\label{forS}
    S(x,t,\xi) = \nu\int_0^t{(pH_p(X,p,\tau) - H(X,p,\tau))d\tau}
\end{equation}
where $X$ and $p$ satisfy 
\begin{equation}\label{withnu}
    \begin{cases}
        \frac{dX}{d\sigma} = \nu H_p(X,p), \quad X|_{\sigma = 0} = \xi \\
        \frac{dp}{d\sigma} = -\nu H_x(X,p), \quad p|_{\sigma = 0} = p_o \in \mathbb{R}
    \end{cases}
\end{equation}
Solving (\ref{withnu}), we have
\begin{equation}
    \begin{cases}
        X = \xi + \nu \sigma H_p(\xi,p_0) + O(\nu ^2) \\
        p = p_0 - \nu \sigma H_x(\xi, p_0)+ O(\nu ^2)
    \end{cases}
\end{equation}
Substituting the expressions for $X$ and $p$ into (\ref{forS}) with $\nu \to +0$, eventually we obtain 
\begin{equation}
\begin{gathered}
    S = \nu t' (pH_p(\xi,p_0) - H(\xi,p_0)) + O(\nu ^2) = \\ = t (pH_p(\xi,p_0) - H(\xi,p_0)) + O(t ^2)
\end{gathered}
\end{equation}

Let $X|_{\sigma = t'} = x \in \mathbb{R}$. Denote $Y = \frac{x-\xi}{\nu t'}$. Expanding $p_0$ into a power series with respect to $\nu$, we have $p_0 = P_0 + \nu P_1 + ...+ \nu^n P_n$. Then $P_0 = P_0(Y, \xi)$ can be found from the relation
\begin{equation}\label{forP0}
    \frac{x-\xi}{\nu t'} = H_p(\xi, P_0)
\end{equation}
In our case $H = x^2p^2$ and
\begin{equation}\label{forS123}
    S = \nu \int_0^{t'} X^2 p^2 d\sigma
\end{equation}
Relation (\ref{forP0}) looks as 
\begin{equation}\label{forPO1}
    \frac{x-\xi}{\nu t'} = 2\xi^2 P_0
\end{equation}
If $|x-\xi| \leq c\nu$, then $P_0$ is a smooth function as long as $\xi \ne 0$. As was mentioned in (\ref{DELTA111}), given $\xi = 0$, the fundamental solution is the Dirac delta function and function $S$ has the form
\begin{equation}
    S(x,\beta) = \frac{x^2}{2(1-\beta)}, \quad \beta \to 1-0
\end{equation}
If $\xi \ne 0$, we have
\begin{equation}
    P_0 = \frac{Y}{2\xi^2}
\end{equation}
Thus, $P_0 \sim \frac{1}{\nu}$ for each fixed $\xi$ and $\nu \to +0$. This implies that we need to set $P = \nu p$ in (\ref{withnu}) in order to obtain $X$ and $P$ by using regular perturbation methods. The system (\ref{withnu}) now takes the form
\begin{equation}\label{withnu1}
    \begin{cases}
        \frac{dX}{d\sigma} = 2X^2P \\
        \frac{dP}{d\sigma} = -2XP^2
    \end{cases}
\end{equation}
with boundary conditions 
\begin{equation}
    X|_{\sigma = 0} = \xi, \quad X|_{\sigma = t'} = x
\end{equation}
According the regular perturbation theory, we seek the solution of such system in the form
\begin{equation}
    X = X_0 + \nu X_1 + ... + \nu^N X_N, \quad P = P_0 + \nu P_1 + ... + \nu^N P_N
\end{equation}
After substituting the expansions for $X$ and $P$ in (\ref{withnu1}), we obtain systems for $X_0$ and $P_0$
\begin{equation}\label{X0P0}
    \begin{cases}
        \frac{dX_0}{d\sigma} = 2X_0^2P_0 \quad X_0|_{\sigma = 0} = \xi\\
        \frac{dP_0}{d\sigma} = -2X_0P_0^2 \quad X_0|_{\sigma = t'} = x
    \end{cases}
\end{equation}
for $X_1$ and $P_1$
\begin{equation}\label{X1P1}
    \begin{cases}
        \frac{dX_1}{d\sigma} = 4X_0 X_1 P_0 + 2 X_0^2 P_1 \quad X_1|_{\sigma = 0} = 0\\
        \frac{dP_1}{d\sigma} = -2X_1P_0^2 - 4X_0 P_0 P_1 \quad X_1|_{\sigma = t'} = 0
    \end{cases}
\end{equation}
and so on. The function $S$ in this case has the form 
\begin{equation}
    S = \nu^{-1}S_0 + S_1 + \nu S_2 + O(\nu^2)
\end{equation}
where $S_j$ can be found from (\ref{forS123}).
\section{Conclusion. Investigation of the Hamilton system in the general one-dimensional case}

$\quad$ The achieved results can be transferred to the case of more general single-dimensional Hamiltonians
\begin{equation}
    H(x,p) = x^2a^2(x)p^2,\quad a \ne 0
\end{equation}
Let us denote $A(x) = xa(x)$. The characteristic system for the corresponding Hamilton-Jacobi equation thus appears the following way
\begin{equation}\label{HAMGEN}
    \begin{cases}
        \dot{x} = H_p = 2A^2p, \quad x|_{t=0} = x_0 \\
        \dot{p} = -H_x = -2AA'p^2, \quad p|_{t=0} = p_0 = \beta(x_0-\xi-y)
    \end{cases}
\end{equation}
When the second equation of the system (\ref{HAMGEN}) is divided by the first one, the ordinary differential equation is obtained
\begin{equation*}
    \frac{dp}{dx} = -\frac{A'}{A}p
\end{equation*}
or the equivalent to it
\begin{equation*}
    \frac{dp}{d} + \frac{dA}{A} = 0
\end{equation*}
Integrating the both parts, we get
\begin{equation}\label{FIRST1}
    p(x_0,t)A(x(x_0,t)) = C_0
\end{equation}
The constant $C_0$ is defined by the initial conditions and after the insertion of $y = -\frac{\beta(x_0-\xi)}{1-\beta}$ has the form
\begin{equation}
    C_0 = p_0 A(x_0) = \frac{\beta x_0a(x_0)(x_0-\xi)}{1-\beta}
\end{equation}
Expressing $p = \frac{C_0}{A}$ from (\ref{FIRST1}) and inserting it in the first equation of the system  (\ref{HAMGEN}),  the ordinary differential equation is obtained
\begin{equation*}
    \frac{dx}{dt} = 2AC_0
\end{equation*}
Separating the variables and integrating, taking into account the initial data, we have
\begin{equation}\label{XXX}
    \int_{x_0}^x \frac{dy}{ya(y)} = \frac{2\beta t x_0a(x_0)(x_0-\xi)}{1-\beta}
\end{equation}
Let us indicate a couple of factors, that could be find out from the expression (\ref{XXX}):

1) $x=0 \Leftrightarrow x_0 = 0$. Inasmuch as with the fixed $t$ and $0<\beta<1$ there is a bounded expression on the right, and the integral on the left diverges at zero value of one of the integration limits, the equality itself (\ref{XXX}) is possible only in the case when the other integration limit is equal to zero.

2) With $\beta \to 1 - 0$, on the contrary, the expression on the right will increase infinitely, and since the integral on the left is bounded, equality is possible only for $x_0 = 0$ or $x_0 = \xi$.

Now let us estimate the properties of the Jacobian $\frac{\partial x}{\partial x_0}$. The integral in the expression (\ref{XXX}) is a function of $ x $ and $ x_0 $. We denote
\begin{equation*}
    g(x,x_0) = \int_{x_0}^x \frac{dx}{xa(x)}
\end{equation*}
Differentiate both sides of the expression (\ref{XXX}) by $ x_0 $
\begin{equation*}
    \frac{\partial g}{\partial x}\frac{\partial x}{\partial x_0} + \frac{\partial g}{\partial x_0} = 2t\frac{\partial C_0}{\partial x_0}
\end{equation*}
Thus, the expression for Jacobian gets the following form
\begin{equation}\label{GENJACOB}
    \frac{\partial x}{\partial x_0} = \frac{1}{\frac{\partial x}{\partial x_0}}\biggm(2t\frac{\partial C_0}{\partial x_0} - \frac{\partial g}{\partial x_0}\biggm)
\end{equation}
According to the rule of differentiation of the integral with the variable limits of integration, we have
\begin{equation*}
    \frac{\partial g}{\partial x} = \frac{1}{xa(x)}, \quad \frac{\partial g}{\partial x_0} = -\frac{1}{x_0a(x_0)}
\end{equation*}
The derivative of a constant is equal to
\begin{equation*}
    \frac{\partial C_0}{\partial x_0} = \frac{\beta}{1-\beta}\big(a(x_0)(2x_0-\xi) + x_0a'(x_0)(x_0-\xi)\big)
\end{equation*}
Substituting the obtained values of partial derivatives in (\ref{GENJACOB}), we have the following
\begin{equation}
        \frac{\partial x}{\partial x_0} = \frac{x a(x)}{x_0 a(x_0)}\biggm(1 + \frac{2\beta t x_0 a(x_0)}{1- \beta}\big(a(x_0)(2x_0-\xi)+x_0a'(x_0)(x_0-\xi)\big)\biggm)
\end{equation}

As it was mentioned before, if $x=0$, than $x_0 = 0$ and backwards. Moreover, according to the condition of the problem $a \ne 0$, which means that the term in front of the parentheses does not affect the equality of the Jacobian to zero. Then it remains to examine the expression in parentheses
\begin{equation}\label{SIGN0}
    1 + \frac{2\beta t x_0 a(x_0)}{1- \beta}\big(a(x_0)(2x_0-\xi)+x_0a'(x_0)(x_0-\xi)\big)
\end{equation} 
The considerations here are similar to the ones that were mentioned before the case $H = x^2p^2$ described above

With $0 < \beta < 1$ the expression (\ref{SIGN0}) also never turns into zero, if we restrict $ x $ and $ \xi $ and consider sufficiently small $ t $.

When $\beta \to 1-0$, then $x_0 = 0$ or $x_0 = \xi$, and for these values the expression (\ref{SIGN0}) will also be different from zero everywhere.

So, the solution at all points will be a smooth function, generally speaking, for sufficiently small $ t $ and in the case of a local consideration of the initial data and the solution, as explained in Section 3.

An obstacle to constructing a solution for finite times is the caustics - the supports of the singularities of the projection of the Lagrangian manifold $ \Lambda $ onto the space $ (x, t) $. A standard way to overcome this obstacle is to use the construction of Maslov tunneling canonical operator \cite{masnaz}. You can also use the construction of a generalized asymptotic solution to the Cauchy problem \cite{dn1}. In this paper, such issues are not considered.

\end{document}